\RequirePackage{fix-cm}
\documentclass[smallextended]{svjour3}       
\smartqed  
\usepackage{graphicx}
\usepackage{amssymb,enumerate}
\usepackage{amsmath}
\usepackage{graphics}
\usepackage{amsfonts}
\usepackage{amsmath, enumerate}
\usepackage{footnote}
\usepackage{multicol}
\usepackage{booktabs}
\usepackage{tikz}
\usepackage[flushleft]{threeparttable}
\usepackage[font=small,labelfont=bf]{caption}
 \usepackage{cite}
\begin{document}

\title{Existence and stability analysis of solution for fractional delay differential equations}

\subtitle{}


\author{Faruk Develi \and Okan Duman}


\institute{F. Develi \at
              Department of Mathematics, Yildiz Technical University, 34210 Istanbul, Turkey \\
              \email{fdeveli@yildiz.edu.tr}           
           \and
           O. Duman \at
              Department of Mathematics, Yildiz Technical University, 34210 Istanbul, Turkey\\
							\email{oduman@yildiz.edu.tr}
}
\date{Received: date / Accepted: date}
\maketitle

\begin{abstract}
In this article, we give some results for fractional-order delay differential equations. In the first result, we prove the existence and uniqueness of solution by using Bielecki norm effectively. In the second result, we consider a constant delay  form of this problem. Then we apply Burton's method to this special form to prove that there is only one solution. Finally, we prove a result regarding the Hyers-Ulam stability of this problem. Moreover, in these results, we omit the conditions for contraction constants seen in many papers.
\end{abstract}

\section{Introduction}
Fractional differential equations appear in various fields. For examples, it’s handled  in engineering with physical processes such as thermodynamics, polymer rheology, and mechanics as well as control theory and technical sciences like biophysics \cite{Ga, G, P, VS, H, T}. The biggest factor in fractional differentiations being more popular than classical ones lately is that it is effective in explaining real-world problems.  But since there is no general technique for  obtaining  solutions of the dynamic systems specified by fractional calculus, existence and uniqueness theorems have an important place in the literature \cite{ABR, Zhou, L, DR, Z, Ragusa, DF, BL}.

The concept of stability for functional equations was first introduced by Ulam at a conference in 1940. After Hyers' first contribution to Ulam's work in 1941, this type of stability concept came to be known as Hyers-Ulam stability. Obloza is the first author to study the this type of stability of linear differential equations \cite{Obloza}. Later, the concept of Hyers-Ulam stability is discussed in many topics such as ordinary differential equations, partial differential equations, and delay differential equations \cite{Jung, LP, Otrocol, Rus, WFZ}. This type stability for fractional-order differential equations with respect to Caputo derivative are investigated by many authors \cite{WLY, WLZ, WZF}. In particular, existence-uniqueness and Hyers-Ulam type stability results regarding delay differential equations were investigated in the papers \cite{Otrocol, WZ, CS3, BP}. In \cite{Otrocol}, results for the delayed form of differential equations in the classical sense were obtained by the Picard operator method, and inspired by this paper, existence-uniqueness and stability results for fractional Caputo derivative were obtained in \cite{WZ}.

In this article, we investigate the existence and uniqueness of
solutions and Hyers-Ulam type stability for the following fractional-order delay differential equation in the sense of Caputo, motivated by \cite{B1, B2, BP, Otrocol, WZ}.
\begin{align} \label{m1}
\begin{cases}
^{c}D^{\alpha}\upsilon(t)=f(t,\upsilon(t),\upsilon(g(t))) \quad &t\in [0,T] \\
~~~~~~ \upsilon(t)=\phi(t) \quad &t\in [-h,0].
\end{cases}
\end{align}
where $f\in C([0,T] \times \mathbb{R}^{2},\mathbb{R})$, $\phi \in C([-h,0],\mathbb{R})$, $g \in C([0,T], [-h,T])$ verifying $g(t) \leq t$ and $^{c}D^{\alpha}$ is the fractional derivative of order $\alpha \in (0,1)$ in the sense of Caputo.

The paper is structured as follows: Section \ref{pre} introduces the general notion of Caputo fractional-order derivative and Hyers-Ulam stability. We state here basic properties of Caputo fractional-order derivative and useful inequalities. In Section \ref{existence}, we investigate the existence and uniqueness of the solution to this problem. In the existence and uniqueness theorems, contractivity constants are one of the most important tools to prove main theory. For this reason, some conditions are put on these constants as a hypothesis before main results are given. In Theorem \ref{main}, we obtain the existence and uniqueness of the solution under the hypothesis that the function on the righthand side of our problem satisfies the Lipschitz condition with respect to the second and third variables. Here, we omit some of the conditions in the article \cite{WZ} and give the proof of the existence and uniqueness theorem by using the Bielecki norm more effectively. In other words, we obtain our result without the need for contraction constants. Later, we give the existence and uniqueness theorem for a special case of the problem (\ref{m1}) by applying the technique named progressive contractions, which is introduced by Burton \cite{B1, B2, BP}. In Burton's progressive contractions, the interval studied is divided into n-equal parts of a certain length. For the first interval, contraction mapping defining through the hypothesis, and a unique solution is obtained. Then, this solution is considered as the initial function and a new contraction mapping is defined for the second interval, and a solution is obtained here. By continuing this process n-steps, we have a unique solution for the whole domain. For more details and other applications of it, see also \cite{R33, Otrocol2} and references therein. This technique allows us to omit the Lipschitz condition stated in Theorem \ref{main} with respect to the third variable.  In Section \ref{stability}, we focus our attention on Hyers-Ulam stability for the problem. Here we obtain the Hyers-Ulam stability result for the equation (\ref{m1}) using Picard operators theory, and then we give an alternative proof without using techniques such as Picard operator theory and Gronwall type inequalities.  Finally, we give examples to illustrate our results in Section \ref{example}.

\section{Preliminaries}\label{pre}
In this section, we present some notations, definitions, and preliminary facts used throughout this paper.
\begin{definition} \cite{P, KST}
The Riemann–Liouville integral of order $\alpha>0$ for the function $\upsilon$ is defined as
\begin{equation*} \label{int}
    I^{\alpha}\upsilon(t)=\frac{1}{\Gamma(\alpha)} \int_{0}^{t}(t-s)^{\alpha-1}\upsilon(s)ds, \quad t\in [0,T],
\end{equation*}
where $\Gamma(\cdot)$ is the Gamma function.
\end{definition}

\begin{definition} \cite{P, KST}
The Caputo derivative of fractional-order $\alpha$ for the function $\upsilon$ is defined as
\begin{equation*}
    D^{\alpha}\upsilon(t)=\frac{1}{\Gamma(n-\alpha)}\int_{0}^{t}(t-s)^{n-\alpha-1}\upsilon^{(n)}(s)ds, \quad t\in [0,T],
\end{equation*}
where $n=[\alpha]+1$ and $[\alpha]$ denotes the integer part of $\alpha$.
\end{definition}
\begin{definition}
The equation (\ref{m1}) is Hyers-Ulam stable if there exists a real number $c>0$ such that for each $\epsilon>0$ and for each solution $\vartheta\in C([-h,T],\mathbb{R})$ to the inequality
\begin{equation}\label{stable}
\left\vert  D^{\alpha}\vartheta(t)-f(t,\vartheta(t),\vartheta(g(t))) \right \vert \leq \epsilon,\quad t \in [0,T], 
\end{equation}
there exists a solution $\upsilon \in C([-h,T],\mathbb{R})$ to the equation (\ref{m1}) with
\begin{equation*}
\big\vert \vartheta(t)-\upsilon(t) \big \vert \leq c\epsilon,\quad t \in [-h,T].
\end{equation*}
\end{definition}

\begin{remark} \label{R1}
A function $\vartheta \in C([0,T],\mathbb{R})$ is a solution of inequality (\ref{stable}) if and only if there exists a function $\Psi \in C([0,T],\mathbb{R})$ such that
\begin{enumerate}[i)]
      \item  $\big\vert \Psi(t) \big \vert \leq \epsilon \quad \text{for all} \quad t\in [0,T]$,
      \item  $^{c}D^{\alpha} \vartheta(t) = f(t,\vartheta(t),\vartheta(g(t)))+\Psi(t) \quad \text{for all} \quad t \in [0,T]$.
    \end{enumerate}
\end{remark}  

\begin{remark} \label{R2}
If $\vartheta \in C([0,T],\mathbb{R})$ is a solution of the inequality (\ref{stable}), then it is a solution to the following integral inequality:
\begin{align*}
\bigg\vert \vartheta(t)-\vartheta(0)-\frac{1}{\Gamma(\alpha)}\int_{0}^{t} (t-s)^{\alpha -1} f(s,\vartheta(s),\vartheta(g(s))ds\bigg\vert \leq \frac{\epsilon T^{\alpha}}{\Gamma(\alpha+1)}
\end{align*}
for all $t \in [0,T]$.
\end{remark}

Now we give the following simple inequality which is useful for our results.
For $\alpha,\tau>0$,
\begin{equation} \label{tool}
    \int_{0}^{t} (t-s)^{\alpha -1} e^{\tau s}ds \leq \frac{e^{\tau t}}{\tau^{\alpha}} \Gamma(\alpha).
\end{equation}
Actually, by substituting $z=t-s$ in the integral expression above, we get 
\begin{align*}
    \int_{0}^{t} (t-s)^{\alpha -1} e^{\tau s}ds &= e^{\tau t} \int_{0}^{t}z^{\alpha-1}e^{-\tau z}dz  \\
    &= \frac{e^{\tau t}}{\tau^{\alpha}} \int_{0}^{\tau t}x^{\alpha-1}e^{-x}dx  \quad \text{substituting $x=\tau z$} \\
    &\leq \frac{e^{\tau t}}{\tau^{\alpha}} \int_{0}^{\infty}x^{\alpha-1}e^{-x}dx=\frac{e^{\tau t}}{\tau^{\alpha}} \Gamma(\alpha).
\end{align*}

\begin{definition} \cite{R1,R2}
Let $(X,d)$ be a metric space. An $\mathcal{A}:X \rightarrow X$ is a Picard operator if there exists $x^*\in X$ such that (i) $F_{\mathcal{A}}=\{x^*\}$ where $F_{\mathcal{A}}=\{x \in X : \mathcal{A}(x)=x\}$ is the fixed point set of $\mathcal{A}$; (ii) the sequence $(\mathcal{A}^n(x_{0}))_{n \in \mathbb{N}}$ converges to $x^*$ for all $x_{0}\in X$.
\end{definition}

\begin{lemma}\cite{R1,R2} \label{AGL}
Let $(X,d,\leq)$ be an ordered metric space and $\mathcal{A}:X \rightarrow X$ be an increasing Picard operator $(F_{\mathcal{A}}=\{x^*\})$. Then, for $x \in X$, $x \leq \mathcal{A}(x)$ implies $x\leq x^*$ while $x \geq \mathcal{A}(x)$ implies $x\geq x^*$.
\end{lemma}

\begin{lemma}\cite{Henry} \label{Gronwall}
Let $\vartheta:[0,T] \rightarrow [0,\infty)$ be a real function and $w$ be nonnegative, locally integrable function on $[0,T]$. If there are constants $k>0$ and $0<\alpha<1$ such that
\begin{equation*}
    \vartheta(t)\leq \omega(t)+k\int_{0}^{t} (t-s)^{-\alpha} \vartheta(s)ds,
\end{equation*}
then there exists a constant $\delta=\delta(\alpha)$ such that 
\begin{equation*}
    \vartheta(t)\leq \omega(t)+\delta k \int_{0}^{t} (t-s)^{-\alpha} \omega(s)ds.
\end{equation*}
\end{lemma}

\section{Existence and Uniqueness results}\label{existence}
In this section, by overcoming the limitations like contractivity constants, we present the existence and uniqueness results of solution for the problem (\ref{m1}).
\begin{theorem} \label{main}
Suppose that 
\begin{enumerate}
    \item[(C1)] $f\in C([0,T]\times \mathbb{R}^{2},\mathbb{R})$, $g\in C([0,T],[-h,T])$ \text{verifying} $g(t)\leq t$ on $[0,T]$.
    \item[(C2)] There is a constant $L>0$ such that
    \begin{equation*}
        \left\vert f(t,\upsilon_{1},\vartheta_{1})-f(t,\upsilon_{2},\vartheta_{2}) \right \vert \leq L \left ( \left\vert \upsilon_{1}-\upsilon_{2}\right\vert +\left\vert \vartheta_{1}-\vartheta_{2} \right\vert \right)
    \end{equation*}
    for all $\upsilon_{i},\vartheta_{i}\in \mathbb{R}$ $(i=1,2)$ and $t\in [0,T]$.
\end{enumerate}
Then the problem (\ref{m1}) has a unique solution.
\end{theorem}

\begin{proof}
We first convert the problem (\ref{m1}) into a fixed point problem. In this sequel, we consider the operator
\begin{equation*}
    \mathcal{F}:C([-h,T],\mathbb{R}) \rightarrow C([-h,T],\mathbb{R})
\end{equation*}
defined by
\begin{equation*}
\mathcal{F}\upsilon(t)=
\begin{cases}
\phi(t), \quad  &t \in [-h,0] \\
\phi(0)+ \int_{0}^{t} \frac{(t-s)^{\alpha -1}}{\Gamma(\alpha)} f(s,\upsilon(s),\upsilon(g(s))ds, \quad  &t \in [0,T].
\end{cases}
\end{equation*}
Then our aim is reduced to finding a unique fixed point of $\mathcal{F}$. Let consider the Banach space $X:=C([-h,T],\mathbb{R})$ endowed with the following Bielecki norm
\begin{equation}\label{Bielecki}
    \left\Vert \upsilon \right\Vert_{B}=\max_{t\in [-h,T]} \left\vert \upsilon(t) \right\vert e^{-\tau t}.
\end{equation}
To achieve our aim, we show that $\mathcal{F}$ is a contraction mapping on $(X,\left\Vert \cdot \right\Vert_{B})$. For all $\upsilon(t),\vartheta(t) \in X$, $\mathcal{F}\upsilon(t)=\mathcal{F}\vartheta(t)$ if $t \in [-h,0]$, then we take $t \in [0,T]$. Hence
\begin{align*}
\big\vert \mathcal{F}\upsilon(t)&-\mathcal{F}\vartheta
(t)\big\vert \\
\leq& \frac{1}{\Gamma(\alpha)}
\int_{0}^{t} (t-s)^{\alpha -1} \big\vert f(s,\upsilon(s),\upsilon(g(s)))-f(s,\vartheta(s),\vartheta(g(s)))\big\vert ds \\
\leq & \frac{L}{\Gamma(\alpha)} \int_{0}^{t} (t-s)^{\alpha -1}e^{\tau s} \bigg( \max_{-h \leq s \leq T}\big\vert \upsilon(s)-\vartheta(s)\big\vert e^{-\tau s} \\
&+ \max_{-h \leq s \leq T}\big\vert \upsilon(g(s))-\vartheta(g(s))\big\vert e^{-\tau s} \bigg) ds \\
\leq & \frac{2L}{\Gamma(\alpha)} \big\Vert \upsilon-\vartheta \big\Vert _{B} \int_{0}^{t} (t-s)^{\alpha -1} e^{\tau s}ds \\
\leq& \frac{2L}{\tau^{\alpha}} \big\Vert \upsilon-\vartheta \big\Vert _{B} e^{\tau t} \quad \text{(by the inequality (\ref{tool}))}.
\end{align*}
Then we obtain that
\begin{equation*}
    \left\Vert \mathcal{F}\upsilon-\mathcal{F}\vartheta \right\Vert_{B} \leq \lambda \left\Vert
\upsilon-\vartheta \right\Vert_{B} \quad \text{where} \quad \lambda = \frac{2L}{\tau^{\alpha}}.
\end{equation*}
If we choose $\tau>0$ large enough so that $\lambda<1$, then there exists a unique fixed point of $\mathcal{F}$ by the Banach Contraction Principle. Thus the proof is complete.
\end{proof}
\begin{remark}
If we take the following special version of the problem (\ref{m1}) by considering $g(t)=t-r$, where $r>0$ is a constant delay then,
\begin{align} \label{m2}
\begin{cases}
^{c}D^{\alpha}\upsilon(t)=f(t,\upsilon(t),\upsilon(t-r)) \quad &t\in [0,T] \\
~~~~~~ \upsilon(t)=\phi(t) \quad &t\in [-r,0].
\end{cases}
\end{align}
Then we prove the existence and uniqueness of solution for the above fractional-order differential equation under the following Lipschitz condition unlike the Lipschitz condition as stated in Theorem \ref{main} by applying progressive contractions.
\end{remark}
Now we state our result as follows.
\begin{theorem} \label{main2}
Let $f: [0,T] \times \mathbb{R}^{2} \rightarrow \mathbb{R}$ be continuous function. Assume that there exist a positive constant $L$ such that
    \begin{equation*}
        \left\vert f(t,\upsilon_{1},\vartheta)-f(t,\upsilon_{2},\vartheta) \right \vert \leq L   \left\vert \upsilon_{1}-\upsilon_{2}\right\vert
    \end{equation*}
    for all $\upsilon_{i},\vartheta \in \mathbb{R}$ $(i=1,2)$ and $t\in [0,T]$.
Then the problem (\ref{m2}) has a unique solution.
\end{theorem}

\begin{proof}
It is obvious that the problem (\ref{m2}) is equivalent to the following integral form:
\begin{equation*}
\upsilon(t)=
\begin{cases}
    \phi(t) \quad &-r \leq t \leq 0 \\
		\phi(0)+ \int_{0}^{t} \frac{(t-s)^{\alpha -1}}{\Gamma(\alpha)}  f(s,\upsilon(s),\upsilon(s-r))ds \quad & 0\leq t\leq T.
\end{cases}
\end{equation*}
To apply progressive contractions, we divide the interval $[0,T]$ into $n$ equal parts which have length $S$ where $0<S<r$ and $nS=T$. That is, the partition is as follows:
\begin{equation*}
   0=S_{0}<S_{1}<\cdots<S_{n}=T,\quad S_{i}-S_{i-1}=S.
\end{equation*}
Also we observe that $ t\leq S_{i+1} \Rightarrow t-r\leq S_{i}$ by the following argument: 
\begin{equation*}
    t\leq S_{i+1} \Rightarrow t-r\leq S_{i+1}-r\leq S_{i+1}-S=S_{i}.
\end{equation*}
\flushleft \textbf{Step 1:} Let $(M_{1}, \left\Vert \cdot \right\Vert_{1})$ be complete normed space of continuous functions $\upsilon:[-r,S_{1}] \rightarrow \mathbb{R}$ with the following norm
\begin{equation*}
    \left\Vert \upsilon \right\Vert_{1}=\max_{t\in [-r,S_{1}]} \left\vert \upsilon(t) \right\vert e^{-\tau t},
\end{equation*}
and we take $\upsilon(t)=\phi(t)$ for $-r \leq t \leq 0$. Define a mapping $\mathcal{F}_{1}:M_{1} \rightarrow M_{1}$ given by
\begin{equation*}
\mathcal{F}_{1}\upsilon(t)=
\begin{cases}
    \phi(t) \quad &-r \leq t \leq 0 \\
		\phi(0)+ \int_{0}^{t} \frac{(t-s)^{\alpha -1}}{\Gamma(\alpha)} f(s,\upsilon(s),\upsilon(s-r))ds \quad & 0\leq t\leq S_{1}.
\end{cases}
\end{equation*}
For $\upsilon(t),\vartheta(t)\in M_{1}$, $\mathcal{F}_{1}\upsilon(t)=\mathcal{F}_{1}\vartheta(t)$ if $t \in [-r,0]$, then we take $t \in [0,S_{1}]$. Hence
\begin{equation*}
    \left\vert \mathcal{F}_{1}\upsilon(t)-\mathcal{F}_{1}\vartheta(t)\right\vert \leq 
\int_{0}^{t} \frac{(t-s)^{\alpha -1}}{\Gamma(\alpha)} \big\vert f(s,\upsilon(s),\upsilon(s-r))-f(s,\vartheta(s),\vartheta(s-r))\big\vert ds.
\end{equation*}
Since $0\leq s\leq S_{1}\Rightarrow (s-r )\in  [-r,0]$ and the definition of $M_{1}$, we have
\begin{align*}
\big\vert \mathcal{F}_{1}\upsilon(t)&-\mathcal{F}_{1}\vartheta(t)\big\vert \\
\leq&\frac{1}{\Gamma(\alpha)} \int_{0}^{t} (t-s)^{\alpha -1} \left\vert f(s,\upsilon(s),\phi(s-r))-f(s,\vartheta(s),\phi(s-r))\right\vert ds \\
\leq&  \frac{L}{\Gamma(\alpha)} \int_{0}^{t} (t-s)^{\alpha -1} e^{\tau s} \bigg( \max_{-h \leq s \leq S_{1}} \big\vert \upsilon(s)-\vartheta(s)\big\vert e^{-\tau s} \bigg) ds \\
\leq& \frac{L}{\Gamma(\alpha)} \left\Vert \upsilon-\vartheta \right\Vert _{1} \int_{0}^{t} (t-s)^{\alpha -1} e^{\tau s}ds \leq \frac{L}{\tau^{\alpha}} \left\Vert \upsilon-\vartheta \right\Vert _{B} e^{\tau t}.
\end{align*}
Consequently, we obtain that
\begin{equation*}
    \left\Vert \mathcal{F}_{1}\upsilon-\mathcal{F}_{1}\vartheta \right\Vert_{1} \leq \lambda \left\Vert
\upsilon-\vartheta \right\Vert_{1} \quad \text{where} \quad \lambda = \frac{L}{\tau^{\alpha}}.
\end{equation*}
By taking $\tau>0$ such that $\lambda<1$, then we have $\mathcal{F}_{1}$ is a contraction mapping and so there exists a unique fixed point $\phi_{1} \in M_{1}$ such that it satisfies the problem (\ref{m2}) on $[-r,S_{1}]$.
\flushleft \textbf{Step 2:} In this step, we extend the interval of Step 1 into $[-r,S_{2}]$. Let $(M_{2}, \left\Vert \cdot \right\Vert_{2})$ be complete normed space of continuous functions $\upsilon:[-r,S_{2}] \rightarrow \mathbb{R}$ with the following norm
\begin{equation*}
    \left\Vert \upsilon \right\Vert_{2}=\max_{t\in [-r,S_{2}]} \left\vert \upsilon(t) \right\vert e^{-\tau t},
\end{equation*}
and we take $\upsilon(t)=\phi_{1}(t)$ for $-r \leq t \leq S_{1}$. Similarly, we define a mapping $\mathcal{F}_{2}:M_{2} \rightarrow M_{2}$ given by
\begin{equation*}
\mathcal{F}_{2}\upsilon(t)=
\begin{cases}
    \phi_{1}(t) \quad &-r \leq t \leq S_{1} \\
		\phi(0)+\int_{0}^{t} \frac{(t-s)^{\alpha -1}}{\Gamma(\alpha)} f(s,\upsilon(s),\upsilon(s-r))ds \quad & S_{1} \leq t\leq S_{2}.
\end{cases}
\end{equation*}
For $\upsilon(t),\vartheta(t)\in M_{2}$, $\mathcal{F}_{2}\upsilon(t)=\mathcal{F}_{2}\vartheta(t)$ if $t \in [-r,S_{1}]$, then we take $t \in [S_{1},S_{2}]$. Thus
\begin{equation*}
    \left\vert \mathcal{F}_{2}\upsilon(t)-\mathcal{F}_{2}\vartheta(t)\right\vert \leq 
\int_{0}^{t} \frac{(t-s)^{\alpha -1}}{\Gamma(\alpha)} \big\vert f(s,\upsilon(s),\upsilon(s-r))-f(s,\vartheta,\vartheta(s-r))\big\vert ds.
\end{equation*}
Note that $0 \leq s\leq S_{2}\Rightarrow (s-r )\in  [-r,S_{1}]$. By considering the definition of $M_{2}$, we may write
\begin{align*}
\big\vert \mathcal{F}_{2}\upsilon(t)&-\mathcal{F}_{2}\vartheta(t)\big\vert \\
\leq&  \frac{1}{\Gamma(\alpha)} \int_{0}^{t} (t-s)^{\alpha -1}  \left\vert f(s,\upsilon(s),\phi_{1}(s-r))-f(s,\vartheta(s),\phi_{1}(s-r))\right\vert ds \\
\leq&  \frac{L}{\Gamma(\alpha)} \int_{0}^{t} (t-s)^{\alpha -1} e^{\tau s} \bigg( \max_{-h \leq s \leq S_{2}} \big\vert \upsilon(s)-\vartheta(s)\big\vert e^{-\tau s} \bigg) ds \\
\leq& \frac{L}{\Gamma(\alpha)} \left\Vert \upsilon-\vartheta \right\Vert _{2} \int_{0}^{t} (t-s)^{\alpha -1} e^{\tau s}ds \leq \frac{L}{\tau^{\alpha}} \left\Vert \upsilon-\vartheta \right\Vert _{B} e^{\tau t}
\end{align*}
and consequently we have
\begin{equation*}
    \left\Vert \mathcal{F}_{2}\upsilon-\mathcal{F}_{2}\vartheta \right\Vert_{2} \leq \lambda \left\Vert
\upsilon-\vartheta \right\Vert_{2},
\end{equation*}
where $\lambda$ is as stated in Step 1. Therefore $\mathcal{F}_{2}$ has a unique fixed point  $\phi_{2}$ in $M_{2}$ such that it satisfies the problem (\ref{m2}) on $[-r,S_{2}]$.
\flushleft \textbf{Step 3:} By continuing this process to $n^{th}$ Step, we can find a continuous mapping $\phi_{n}$ as in the other Steps, which is the unique solution for the problem (\ref{m2}) on $[-r,S_{n}]=[-r,T]$.
\end{proof}

\section{Hyers-Ulam stability result}
In this section we give a result on the Hyers-Ulam stability of the first equation in the problem (\ref{m1}).

\begin{theorem} \label{stability}
Assume that the conditions (C1) and (C2) are fulfilled. Then the first equation of the problem (\ref{m1}) is Hyers-Ulam stable.
\end{theorem}

\begin{proof}
Let $\vartheta$ be a solution to (\ref{stable}). We indicate $\upsilon$ as a unique solution to the following problem by Theorem \ref{main},
\begin{equation*} 
\begin{cases}
^{c}D^{\alpha}\upsilon(t)=f(t,\upsilon(t),\upsilon(g(t))) \quad &t\in [0,T] \\
~~~~~~\upsilon(t)=\vartheta(t) \quad &t\in [-h,0]
\end{cases}
\end{equation*}
It follows we have
\begin{equation*}
\upsilon(t)=
\begin{cases}
\vartheta(t), \quad  &t \in [-h,0] \\
\vartheta(0)+ \int_{0}^{t} \frac{(t-s)^{\alpha -1}}{\Gamma(\alpha)} f(s,\upsilon(s),\upsilon(g(s))ds, \quad  &t \in [0,T].
\end{cases}
\end{equation*}
Obviously, we also have from Remark \ref{R2}
\begin{equation*}
    \Big\vert \vartheta(t)-\vartheta(0) -\frac{1}{\Gamma(\alpha)}\int_{0}^{t} (t-s)^{\alpha -1} f(s,\vartheta(s),\vartheta(g(s))ds\Big\vert \leq \frac{\epsilon T^{\alpha}}{\Gamma(\alpha+1)}
\end{equation*}
for all $t\in[0,T]$, and $\big\vert \vartheta(t)-\upsilon(t) \big\vert=0$ for all $t\in[-h,0]$. For all $t\in[0,T]$, we obtain from the condition (C2) that
\begin{align}
    \big\vert \vartheta&(t)-\upsilon(t) \big\vert \nonumber\\ \leq& \Big\vert \vartheta(t)-\vartheta(0) -\frac{1}{\Gamma(\alpha)}\int_{0}^{t} (t-s)^{\alpha -1} f(s,\vartheta(s),\vartheta(g(s))ds\Big\vert \nonumber\\
    &+ \frac{1}{\Gamma(\alpha)}\int_{0}^{t} (t-s)^{\alpha -1} \big\vert f(s,\vartheta(s),\vartheta(g(s))-f(s,\upsilon(s),\upsilon(g(s)) \big\vert \nonumber\\
    \leq& \frac{\epsilon T^{\alpha}}{\Gamma(\alpha+1)} + \frac{L}{\Gamma(\alpha)}\int_{0}^{t} (t-s)^{\alpha -1} \Big( \big\vert \vartheta(s)-\upsilon(s) \big\vert +\big\vert \vartheta(g(s))-\upsilon(g(s)) \big\vert \Big). \label{eq}
\end{align}
For $z\in C([-h,T],\mathbb{R}^+)$, we define the operator
\begin{equation*}
     \mathcal{A}:C([-h,T],\mathbb{R}^+) \rightarrow C([-h,T],\mathbb{R}^+)
\end{equation*}
given by 
\begin{equation*}
\mathcal{A}(z)(t)=
    \begin{cases}
        0 &\quad t \in [-h,0] \\
       \frac{\epsilon T^{\alpha}}{\Gamma(\alpha+1)} + \frac{L}{\Gamma(\alpha)}\int_{0}^{t} (t-s)^{\alpha -1} \Big( z(s)+z(g(s)) \Big)ds &\quad t\in [0,T].
    \end{cases}
\end{equation*}
To prove that $\mathcal{A}$ is a Picard operator, we show that $\mathcal{A}$ is a contraction mapping with the Bielecki norm given in (\ref{Bielecki}). For $z,\tilde{z}\in C([-h,T],\mathbb{R}^+)$, we have 
\begin{align*}
  \big\vert \mathcal{A}z-\mathcal{A}\tilde{z} \big\vert &\leq \frac{L}{\Gamma(\alpha)}\int_{0}^{t} (t-s)^{\alpha -1} \Big( \big\vert z(s)-\tilde{z}(s) \big\vert + \big\vert z(g(s))-\tilde{z}(g(s)) \big\vert \Big)ds \\
  &\leq \frac{2L}{\Gamma(\alpha)} \big\Vert z-\tilde{z} \big\Vert _{B} \int_{0}^{t} (t-s)^{\alpha -1} e^{\tau s}ds \leq \frac{2L}{\tau^{\alpha}} \big\Vert z-\tilde{z}\big\Vert _{B} e^{\tau t}
\end{align*}
which implies that 
\begin{equation*}
    \big\Vert \mathcal{A}z-\mathcal{A}\tilde{z} \big\Vert_{B} \leq \lambda \big\Vert
z-\tilde{z} \big\Vert_{B} \quad \text{where} \quad \lambda = \frac{2L}{\tau^{\alpha}}.
\end{equation*}
Choosing an appropriate real number $\tau>0$ such that  $\lambda<1$, we get that $\mathcal{A}$ is a contraction mapping with respect to the Bielecki norm $\left\Vert \cdot \right\Vert_{B}$ on $C([-h,T],\mathbb{R}^+)$. Hence $\mathcal{A}$ is a Picard operator such that $F_{\mathcal{A}}=\{z^*\}$ and the following equality holds by the Banach contraction principle  
\begin{align*}
    z^{*}(t)= \frac{\epsilon T^{\alpha}}{\Gamma(\alpha+1)} + \frac{L}{\Gamma(\alpha)}\int_{0}^{t} (t-s)^{\alpha -1} \Big( z^{*}(s)+z^{*}(g(s)) \Big)ds  
\end{align*}
for $t\in [0,T]$. To show that $z^{*}$ is increasing, we denote $m:=\min_{t\in [0,T]}[z^{*}(t)+z^*(g(t))]\in \mathbb{R}^{+}$. For $0\leq t_{1}<t_{2} \leq T$, we have
\begin{align*}
    z^{*}(t_{2})-z^{*}(t_{1})=&\frac{L}{\Gamma(\alpha)}\int_{0}^{t_{1}} \Big( (t_{2}-s)^{\alpha -1}-(t_{1}-s)^{\alpha -1} \Big) \Big( z^{*}(s)+z^{*}(g(s)) \Big)ds \\
    &+ \frac{L}{\Gamma(\alpha)}\int_{t_{1}}^{t_{2}} (t_{2}-s)^{\alpha -1} \Big( z^{*}(s)+z^{*}(g(s)) \Big)ds \\
    \geq&  \frac{mL}{\Gamma(\alpha)}\int_{0}^{t_{1}}  \Big( (t_{2}-s)^{\alpha -1}-(t_{1}-s)^{\alpha -1} \Big) ds \\
    &+ \frac{mL}{\Gamma(\alpha)}\int_{t_{1}}^{t_{2}} (t_{2}-s)^{\alpha -1} ds \\
    =& \frac{mL}{\Gamma(\alpha+1)} (t_{2}^{\alpha}-t_{1}^{\alpha}) > 0.
\end{align*}
Then we can say that the solution $z^*$ is increasing and so $z^*(g(t)) \leq  z^*(t)$ due to $g(t)\leq t$. It follows that
\begin{equation*}
    z^*(t) \leq \frac{\epsilon T^{\alpha}}{\Gamma(\alpha+1)} + \frac{2L}{\Gamma(\alpha)}\int_{0}^{t} (t-s)^{\alpha -1}  z^{*}(s) ds.
\end{equation*}
Applying Lemma \ref{Gronwall} to the above inequality, we obtain that 
\begin{equation*}
     z^*(t) \leq \frac{\epsilon T^{\alpha}}{\Gamma(\alpha+1)} \bigg( 1+ \frac{2\delta LT^{\alpha}}{\Gamma(\alpha+1)}\bigg)
\end{equation*}
for all $t\in [-h,T]$. In particular, if we choose $z=\big\vert \vartheta-\upsilon \big\vert$ in (\ref{eq}), then $z \leq \mathcal{A}z$. So, we have $z \leq z^*$ by Lemma \ref{AGL} since $\mathcal{A}$ is an increasing Picard operator. Consequently, we have
\begin{equation*}
    \big\vert \vartheta(t)-\upsilon(t) \big\vert \leq c\epsilon \quad \text{where} \quad c=\frac{ T^{\alpha}}{\Gamma(\alpha+1)} \bigg( 1+ \frac{2\delta LT^{\alpha}}{\Gamma(\alpha+1)}\bigg).
\end{equation*}
Thus the first equation of (\ref{m1}) is Hyers-Ulam stable.
\end{proof}

\begin{proof}[Alternative proof]
By considering the inequality of (\ref{eq}), we have
\begin{align*}
     \big\vert \vartheta(t)-\upsilon(t) \big\vert 
     \leq& \frac{\epsilon T^{\alpha}}{\Gamma(\alpha+1)} + \frac{L}{\Gamma(\alpha)}\int_{0}^{t} (t-s)^{\alpha -1} \Big( \big\vert \vartheta(s)-\upsilon(s) \big\vert \\
     &+\big\vert \vartheta(g(s))-\upsilon(g(s)) \big\vert \Big)ds \\
     \leq & \frac{\epsilon T^{\alpha}}{\Gamma(\alpha+1)} + \frac{2L}{\Gamma(\alpha)} \left\Vert \vartheta-\upsilon \right\Vert _{B} \int_{0}^{t} (t-s)^{\alpha -1} e^{\tau s}ds  \\
     \leq & \frac{\epsilon T^{\alpha}}{\Gamma(\alpha+1)}+\frac{2L}{\tau^{\alpha}} \left\Vert \vartheta-\upsilon \right\Vert _{B} e^{\tau t}.
\end{align*}
Then we have
\begin{equation*}
(1-\lambda)\big\Vert \vartheta-\upsilon \big\Vert_{B} \leq \frac{\epsilon T^{\alpha}e^{\tau h}}{\Gamma(\alpha+1)} \quad \text{where} \quad     \lambda=\frac{2L}{\tau^{\alpha}}.
\end{equation*}
Choosing large enough $\tau>0$ such that $\lambda<1$, we get 
\begin{equation*}
  \big\vert \vartheta(t)-\upsilon(t) \big\vert e^{-\tau t}  \leq \big\Vert \vartheta-\upsilon \big\Vert_{B} \leq  \frac{\epsilon T^{\alpha}e^{\tau  h}}{(1-\lambda) \Gamma(\alpha+1)} 
\end{equation*}
Consequently, we obtain that
\begin{equation*}
        \big\vert \vartheta(t)-\upsilon(t) \big\vert \leq  c \epsilon, \quad  c:=\frac{T^{\alpha}e^{(h+T)\tau}}{(1-\lambda) \Gamma(\alpha+1)}
\end{equation*}
for all $t\in [-h,T]$. Thus the first equation of (\ref{m1}) is Hyers-Ulam stable.
\end{proof}

\section{Examples}\label{example}
\begin{example}
Consider the following fractional-order differential equation
\begin{equation} \label{ex1}
\begin{cases}
^{c}D^{\frac{1}{2}}\upsilon(t)=\frac{\big\vert \upsilon(t)\big\vert}{1+\big\vert \upsilon(t)\big\vert}+\cos{\upsilon(t^{2})} \quad &t\in[0,1] \\
 ~~~~~~\upsilon(t)=t \quad &t\in [-1,0].
\end{cases}
\end{equation}
Let $f(t,\upsilon,\vartheta)=\frac{\big\vert \upsilon \big\vert}{1+\big\vert \upsilon \big\vert}+\cos{\vartheta}$ and $g(t)=t^{2}$. It is clear that
$$\big\vert f(t,\upsilon_{1},\vartheta_{1})-f(t,\upsilon_{2},\vartheta_{2},) \big \vert \leq    \big\vert \upsilon_{1}-\upsilon_{2}\big\vert +\big\vert \vartheta_{1}-\vartheta_{2} \big\vert$$
for all $\upsilon_{i},\vartheta_{i}\in \mathbb{R}$ $(i=1,2)$ and $t\in [0,1]$. Then we obtain from Theorem \ref{main} that the above problem (\ref{ex1}) has a unique solution. In addition, we also obtain that the first equation in this problem is Hyers-Ulam stable by Theorem \ref{stability}.
\end{example}

\begin{example}
Consider the following fractional-order differential equation with a constant delay;
\begin{equation} \label{exp2}
\begin{cases}
^{c}D^{\frac{1}{2}}\upsilon(t)=\sin{\upsilon(t)}+\upsilon^{2}(t-1) \quad &t\in[0,10] \\
 ~~~~~~\upsilon(t)=e^{t} \quad &t\in [-1,0].
\end{cases}
\end{equation}
Let $f(t,\upsilon,\vartheta)=\sin{\upsilon}+\vartheta^{2}$ and $g(t)=t-1$. It is obvious that
$$\left\vert f(t,\upsilon_{1},\vartheta)-f(t,\upsilon_{2},\vartheta) \right \vert \leq    \left\vert \upsilon_{1}-\upsilon_{2}\right\vert$$
for all $\upsilon_{i},\vartheta \in \mathbb{R}$ $(i=1,2)$ and $t\in [0,10]$. Then we obtain from Theorem \ref{main2} that the above problem (\ref{exp2}) has a unique solution without checking the Lipschitz condition with respect to third variable.
\end{example}
\section{Conclusion}
In this article, we consider the main problem in the article \cite{WZ} and motivated by \cite{B1, B2, BP, Otrocol}. We have observed that the contractivity constants appear as an important hypothesis in the main results in \cite{WZ} and in many papers. We omit these hypotheses by using the Bielecki norm more effectively and then show that there is a unique solution. Then we take a special case of our main problem with a constant delay. Here we use Burton's method to show the existence and uniqueness of solution reducing the Lipschitz condition with respect to the third variable. Then we investigate Hyers-Ulam stability of our problem using Picard operators theory, and then we show that the problem can be proven to be stable without the need for techniques such as Picard operator theory and Gronwall-type inequalities. Finally, we give examples to illustrate our results. As can be seen from these results and examples, there is no need for contractivity constant to be less than 1.


\end{document}